\documentclass{article}
\usepackage{amsmath}
\usepackage{graphicx,psfrag,epsf}
\usepackage{enumerate}
\usepackage{natbib}
\usepackage{url} 

\usepackage{comment}

\usepackage[table]{xcolor}
\newcommand{\blind}{0}

\usepackage{soul}
\usepackage{xcolor}     
\usepackage[normalem]{ulem}

\usepackage{calc}
\newsavebox\CBox
\newcommand\hcancel[2][0.5pt]{%
  \ifmmode\sbox\CBox{$#2$}\else\sbox\CBox{#2}\fi%
  \makebox[0pt][l]{\usebox\CBox}%
  \rule[0.5\ht\CBox-#1/2]{\wd\CBox}{#1}}

\begin{document}

\bibliographystyle{agsm}

\def\spacingset#1{\renewcommand{\baselinestretch}%
{#1}\small\normalsize} \spacingset{1}


\if0\blind
{
  \title{\bf Bayes factors and posterior estimation: Two sides of the very same coin}
  \author{Harlan Campbell and Paul Gustafson \\
  Department of Statistics, University of British Columbia
  }
  \maketitle
} \fi

\if1\blind
{
  \bigskip
  \bigskip
  \bigskip
  \begin{center}
    {\LARGE\bf Title}
\end{center}
  \medskip
} \fi

\bigskip

\spacingset{1.45} 

\abstract{Recently, several researchers have claimed that conclusions obtained from a Bayes factor (or the posterior odds) may contradict those obtained from Bayesian posterior estimation.  In this paper, we wish to point out that no such \textcolor{black}{``contradiction''} exists if one is willing to consistently define one's priors and posteriors.  The key for congruence is that the (implied) prior model odds used for testing are the same as those used for estimation.  Our recommendation is simple: If one reports a Bayes factor comparing two models, then one should also report posterior estimates which appropriately acknowledge the uncertainty with regards to which of the two models is correct.
}{\it Keywords:}  Bayes Factor, Bayesian estimation, model averaging.

\section{Introduction}

Recently, several researchers have claimed that conclusions obtained from a Bayes factor (or the posterior odds) may contradict those obtained from Bayesian posterior estimation.  For example, \citet{rouder2018bayesian} discuss what they see as ``two popular Bayesian approaches that may seem incompatible inasmuch as they provide different answers to what appears to be the same question.''  The two approaches in question are referred to as the ``estimation approach'' and the ``Bayes factor approach''.  \citet{wagenmakers2020support} also discuss these two approaches and ask why they result in a ``paradoxical state of affairs''.   \citet{kelter2022evidence} lament how  ``Bayesian interval estimates and hypothesis tests can yield contradictory conclusions'' and \citet{tendeiro2019review} examine an apparent ``mismatch between results from tests and credible intervals.'' \textcolor{black}{\citet{wagenmakers2020principle} go so far as to suggest that, since credible intervals and the Bayes factor are at odds, ``the practice of rejecting [the null] whenever a 95\% interval does not include the null value [...] is not principled and may be misleading in practice.''}


In this paper, we wish to point out that no such ``conflict'', ``paradox'', or ``mismatch'' exists if one is willing to consistently define one's priors and posteriors.  
\textcolor{black}{Specifically, we show that if the same (implied) prior model odds are specified, the Bayes factor approach and the estimation approach are in fact, entirely congruent.}  

Let $\theta$ be the parameter of interest and suppose there are two different models, Model 0 ($M_{0}$) and Model 1 ($M_{1}$), which are \textit{a priori} probable with probabilities $\textrm{Pr}(M_{0})$ and $\textrm{Pr}(M_{1})$, respectively, such that $\textrm{Pr}(M_{0}) + \textrm{Pr}(M_{1}) = 1$.  For each of these models, there is a distinct prior distribution for $\theta$.  For model $i$, let $\pi_{i}(\theta)$ be the prior density for $\theta$, and let $\pi_{i}(\theta|data)$ be the corresponding posterior density such that:
\begin{equation}
\pi_{i}(\theta|data) = \frac{\pi_{i}(\theta)\textrm{Pr}(data|\theta)}{\textrm{Pr}(data|M_{i})},
\nonumber
\end{equation}
\noindent for $i=0,1$, where $\textrm{Pr}(data|\theta)$ is the model distribution of the data given $\theta$, and:
\begin{equation}
\textrm{Pr}(data|M_{i}) = \int{\textrm{Pr}(data|\theta)\pi_{i}(\theta)\textrm{d}\theta}.
\label{eq:pr_data_m}
\end{equation}
\noindent Based on the posterior, one may calculate point estimates for $\theta$, such as the posterior mean or the posterior median, and credible intervals.

One may also be interested in calculating the Bayes factor, $\textrm{BF}_{10}$, which is defined as the ratio of the posterior odds to the prior odds, and can also be defined as the ratio of the marginal likelihoods of the observed data for the two models:
\begin{align}
    \textrm{BF}_{10} &= \frac{\textrm{Pr}(M_{1}|data)}{\textrm{Pr}(M_{0}|data)} / \frac{\textrm{Pr}(M_{1})}{\textrm{Pr}(M_{0})}  \nonumber \\
     &= \frac{\textrm{Pr}(data|M_{1})}{\textrm{Pr}(data|M_{0})}  .
     \nonumber
    \end{align}

\noindent In order to determine which of the two models is more likely to be the true data-generating mechanism, one may consider the posterior model probabilities, $\textrm{Pr}(M_{0}|data)$ and $\textrm{Pr}(M_{1}|data)$:
\begin{align}
\textrm{Pr}(M_{0}|data) &= \frac{\textrm{Pr}(data|M_{0})\textrm{Pr}(M_{0})}{\textrm{Pr}(data|M_{0})\textrm{Pr}(M_{0})+\textrm{Pr}(data|M_{1})\textrm{Pr}(M_{1})} \nonumber \\
&= \frac{\textrm{Pr}(M_{0})}{\textrm{Pr}(M_{0}) + \textrm{BF}_{10} \times \textrm{Pr}(M_{1})}, \quad 
\label{eq:pr_m0_data}
\end{align}
{and}
\begin{align}
\textrm{Pr}(M_{1}|data) &= \frac{\textrm{Pr}(data|M_{1})\textrm{Pr}(M_{1})}{\textrm{Pr}(data|M_{0})\textrm{Pr}(M_{0})+\textrm{Pr}(data|M_{1})\textrm{Pr}(M_{1})} \nonumber \\
&=\frac{\textrm{BF}_{10} \times \textrm{Pr}(M_{1}) }{    \textrm{Pr}(M_{0})+\textrm{BF}_{10} \times \textrm{Pr}(M_{1})},
\label{eq:pr_m1_data}
\end{align}
\noindent as well as the posterior odds, $\textrm{PO}_{10}$, which represent the relative evidence for Model 1 versus Model 0:
\begin{align}
\textrm{PO}_{10} &= \frac{\textrm{Pr}(M_{1}|data)}{\textrm{Pr}(M_{0}|data)} = \frac{\textrm{Pr}(M_{1})}{\textrm{Pr}(M_{0})} \times     \textrm{BF}_{10}.
\label{eq:po_10}
\end{align}

The Bayes factor, $\textrm{BF}_{10}$, can be calculated without defining prior model odds (i.e., without specifying ${\textrm{Pr}(M_{1})}$ and ${\textrm{Pr}(M_{0})}$).  However, in order to determine which model is more likely to be the ``true model'', the Bayes factor must be combined with the prior model odds to obtain the posterior odds.  Despite this fact, the practice of explicitly specifying prior model odds is ``often ignored'' by researchers \citep{tendeiro2019review} who simply quote the Bayes factor as ``the weight of
evidence from the data'' in favour of one model relative to another \citep{o2004kendall}.  \citet{lavine1999bayes} explain in detail why ``such informal use of Bayes factors suffers a certain logical flaw'' that can result in incoherent decisions.

Should a researcher make decisions about which model they believe is most likely to be true citing only the Bayes factor, one must work backwards in order to determine their ``implied prior model odds.''  For instance, if someone believes that $M_{1}$ is more likely to be the \textit{true} model whenever $\textrm{BF}_{10}>1$, and believes that $M_{0}$ is more likely to be the \textit{true} model whenever $\textrm{BF}_{10}<1$, then it follows that such a person has assumed prior model odds of 1:1.  Someone more skeptical of $M_{1}$ might only believe that $M_{1}$ is most likely whenever $\textrm{BF}_{10}>9$, and that $M_{0}$ is most likely whenever $\textrm{BF}_{10}<9$. Such beliefs would correspond to implied prior model odds of 1:9 (i.e., they would necessarily imply that $\textrm{Pr}(M_{0})=0.9$ and $\textrm{Pr}(M_{1})=0.1$)).  
Note that researchers will typically only be willing to conclude with some certainty that $M_{0}$ is the true model if $\textrm{BF}_{10}$ falls bellow some threshold (e.g., if $\textrm{BF}_{10}<1/3$) and only be willing to conclude with some certainty that $M_{1}$ is the true model if $\textrm{BF}_{10}$ falls above some threshold (e.g., if $\textrm{BF}_{10}>3$).  These thresholds are chosen based on both the (implied) prior model odds and on the (implied) relative costs of making a false positive conclusion or false negative conclusion versus remaining indecisive; see \citet{lavine1999bayes}.

Curiously, some researchers adopt different (implied) prior model odds for testing and for estimation.  For instance, based on the idea that there is a fundamental distinction between testing (`is the effect, $\theta$, present or absent?') and estimation (`how big is the effect, $\theta$, assuming it is present?') \citep{wagenmakers2018bayesian}, researchers often use a Bayes factor comparing $M_{1}$ (`effect is present') to a point null $M_{0}$ (`effect is absent') for testing with (implied) prior model odds of 1:1, but then assume $M_{1}$ (`effect is present') is the true model for estimation (with implied prior model odds of 1:0).  As we shall see in the next section with a simple example, this curious practice of adopting different (implied) prior model odds for testing and for estimation is the root cause of the so-called ``paradoxical state of affairs.''

\section{Flipping a possibly biased coin}

As a concrete example, consider observing $X$ ``heads'' out of $N$ coin flips from a possibly biased coin.  This is the same example as considered by \citet{wagenmakers2020support}; see also \citet{puga2015bayes, puga2015bayesian}.  \textcolor{black}{Curiously, while the biased coin example has long been part of statistical folklore, such a coin does not actually exist in the physical world; see \citet{gelman2002you}.}

We assume that the observed coin flip data are the result of a Binomial distribution where the $\theta$ parameter corresponds to the probability of obtaining a ``heads'' such that:
\begin{align}
\textrm{Pr}(data|\theta) =    f_{Binom}(X, N| \theta), \nonumber
\end{align}
\noindent where $f_{Binom}()$ is the Binomial probability mass function.

\begin{figure}
    \centering
    \includegraphics[width=6.2in]{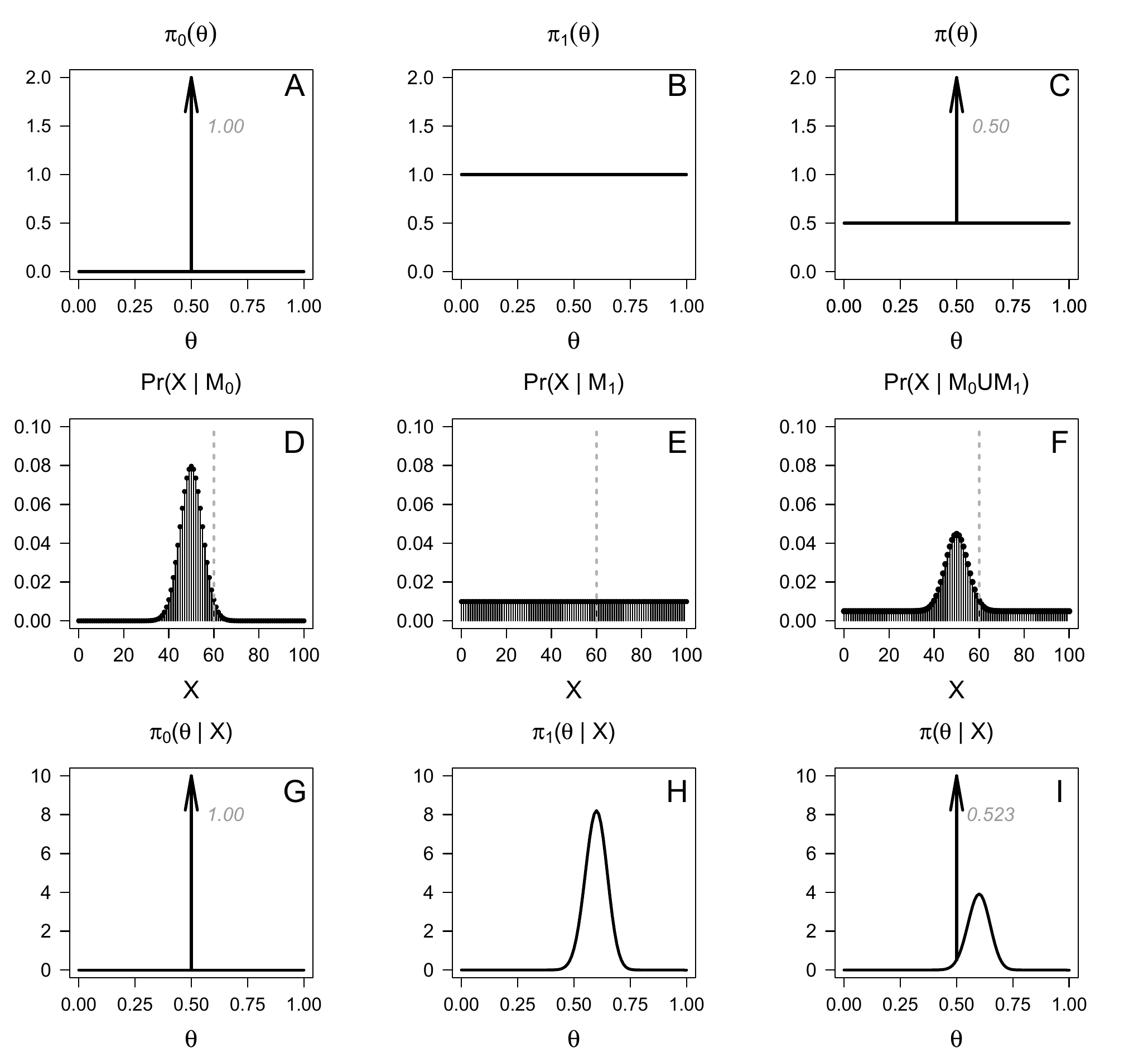}
    \caption{Flipping a possibly biased coin ($X=60$ and $N=100$). The first row shows prior density functions (with the grey numbers next to the arrows corresponding to point mass of the ``spike''); middle row shows probability model for the data (with the grey vertical dashed lines corresponding to observed data of $X=60$); and the lower row shows posterior density functions (with the grey numbers next to the arrows corresponding to point mass of the ``spike'').  Left column corresponds to $M_{0}$, middle column corresponds to $M_{1}$, and right column corresponds to the ``mixed/averaged'' model.}
    \label{fig:3by3}
\end{figure}

 We define two different priors, one for $M_{0}$, and another for $M_{1}$.  For $M_{0}$, a ``point null'' prior states that the coin is fair (i.e., ``heads'' and ``tails'' are equally likely), such that $\theta=0.5$.  For $M_{1}$, the prior states that the probability of a ``heads'' could be anywhere between 0 and 1, with equal likelihood, as described by a Beta distribution:  $\theta \sim \textrm{Beta}(\alpha,\beta)$, where $\alpha=\beta=1$.  The $\textrm{Beta}(1,1)$ distribution is equivalent to a Uniform(0,1) distribution.  See panels A and B in Figure \ref{fig:3by3}. 

The prior density functions are therefore defined as follows:
\begin{equation}
    \pi_{0}(\theta)= {\delta}_{0.5}(\theta) \textrm{  ,} \nonumber
\label{eq:pointnull}
\end{equation}
and
\begin{align}
 \pi_{1}(\theta)=
\begin{cases}
    1 ,& \text{if } \theta\in [0,1] \\
    0,              & \text{otherwise,}
\end{cases}
\label{eq:indicator}
\end{align}
\noindent where $\delta_{0.5}()$ is the Dirac delta function at 0.5 which can be informally thought of as a probability density function which is zero everywhere except at 0.5, where it is infinite.

Suppose we observe $X=60$ ``heads'' out of $N=100$ coin flips.  Then we can easily calculate the following from equation (\ref{eq:pr_data_m}):
\begin{equation}
    \textrm{Pr}(data|M_{0})= f_{Binom}(60, 100, 0.5)=0.0108,  \nonumber
\end{equation}
and
\begin{equation}
    \textrm{Pr}(data|M_{1})= \int_{t=0}^{t=1}f_{Binom}(60, 100, t)dt=0.0099. \nonumber
\end{equation}
\noindent These functions are plotted in panels D and E of Figure \ref{fig:3by3} with the grey vertical dashed line corresponding to observed data of $X=60$.  The ratio of these two numbers is equal to the Bayes factor:
\begin{align}    
   \textrm{BF}_{10}    &= \frac{\textrm{Pr}(data|M_{1})}{\textrm{Pr}(data|M_{0})}  
   = \frac{0.0099}{0.0108} =  0.913.
   \label{eq:BF_coin}
\end{align}
In general, the Bayes factor for this scenario can be computed as:
\begin{align} 
     \textrm{BF}_{10}    &= (N+1)  \binom{N}{X} \theta_{0}^{X}(1-\theta_{0})^{N-X}, \nonumber
\end{align}
where $\theta_{0}=0.5$.  Now suppose the prior probability of each of the two models is equal, such that $\textrm{Pr}(M_{0})=\textrm{Pr}(M_{1}) = 0.5$.  Then, following equations (\ref{eq:pr_m0_data}), (\ref{eq:pr_m1_data}) and (\ref{eq:po_10}), we obtain:
\begin{align}
\textrm{Pr}(M_{0}|data)= \frac{0.5}{0.5 + 0.913 \times 0.5} = 0.523, \nonumber
\end{align}
\begin{align}
\textrm{Pr}(M_{1}|data)= \frac{0.913 \times 0.5 }{    0.5+0.913 \times 0.5} = 0.477, \quad \nonumber
\end{align}
and
\begin{align}
\textrm{PO}_{10} &=  \frac{0.5}{0.5} \times     0.913 = 0.913,
\label{eq:PO_coin}
\end{align}

\noindent which indicates modest support for $M_{0}$ relative to $M_{1}$.

Each of the two models has a corresponding posterior distribution plotted in panels G and H of Figure \ref{fig:3by3}.  If one assumes that $M_{1}$ is the correct model, then the posterior distribution of $\theta$ can be derived analytically (since the Binomial and Beta are conjugate distributions) as:
\begin{align}
    \pi_{1}(\theta|data) &= f_{\textrm{Beta}}(\theta, X+\alpha,N-X+\beta) \nonumber \\
           &= f_{\textrm{Beta}}(\theta, 60+1,100-60+1) \nonumber\\
                      &= f_{\textrm{Beta}}(\theta, 61,41),
\end{align}
\label{eq:post_m1}
\noindent where $f_{\textrm{Beta}}()$ is the Beta probability density function.  The posterior mean, $\hat{\theta}_{1}$, posterior median, $\tilde{\theta}_{1}$, and a 95\% equal-tailed credible interval, $95\%\textrm{CrI}(\theta)_{1}$, can then be calculated as: $\hat{\theta}_{1}= 61/(61+41) = 0.598$, $\tilde{\theta}_{1}= 0.599$, and $95\%\textrm{CrI}(\theta)_{1} = [0.502, 0.691]$.  This credible interval, notably, does not include 0.5.  

If, alternatively,  one assumes that $M_{0}$ is the correct model, then, consequently, the posterior distribution is $\pi_{0}(\theta|data)=\delta_{0.5}(\theta)$.  The posterior mean, median, and 95\% equal-tailed credible interval will be: $\hat{\theta}_{0}= 0.500$, $\tilde{\theta}_{0}= 0.500$, and $95\%\textrm{CrI}(\theta)_{0} = [0.500, 0.500]$.   This credible interval \textit{only} includes 0.500, indicating total certainty with regards to the true value of $\theta$.  \textcolor{black}{As such, it is a very \textit{conservative} 95\% credible interval in the sense that 100\% (instead of exactly 95\%) of the posterior weight is within the interval.}

This coin flip scenario is a good example for illustrating the so-called ``incompatibility'' discussed by \citet{rouder2018bayesian} and others: The $ \textrm{BF}_{10}=0.913$ indicates that the data are more likely to occur under the null, while $95\%\textrm{CrI}(\theta)_{1} = [0.502, 0.691]$ excludes 0.500.  However, there is no reason to expect that the $\textrm{BF}_{10}$ and $95\%\textrm{CrI}(\theta)_{1}$ will be congruent.  If someone decides that $M_{0}$ (`the coin is fair') is more likely because $\textrm{BF}_{10}<1$, necessarily their implied prior model odds must be 1:1.  However, if they also claim that $\text{Pr}(\theta \in 95\%\textrm{CrI}(\theta)_{1}|data) = 0.95$, then they must also believe that $\textrm{Pr}(M_{1})=1$ and $\textrm{Pr}(M_{0})=0$ which is clearly incompatible with the implications made for testing.  To be clear, testing and estimation will only be congruent if the (implied) prior model odds adopted for both are the same.  In the next section, we review two different ways this can be achieved.


\section{Two equivalent approaches}

If someone assumes that $M_{1}$ is the correct model (and therefore reports $95\%\textrm{CrI}(\theta)_{1}$), then there is really no sense for them to compute a Bayes factor (or the posterior odds), since the (implied) prior odds are definitive with respect to which of the two models is ``correct.''  \citet{van2021advantages} explain as follows: If one assumes that $M_{1}$ is the correct model, then the null model is ``deemed false from the outset, and hence no amount of data can either support or undercut it''.

Alternatively, if someone is uncertain about which of the two models is correct, estimates and credible intervals for $\theta$ should take this uncertainty into account.  This can be achieved by either (1) Bayesian model averaging (BMA), or by (2) defining a single ``mixture'' prior.  Both approaches, as we will demonstrate, will deliver the very same results.   

Under the BMA approach, one averages the posterior distributions under each of the two models, $\pi_{0}(\theta|data) $ and   $\pi_{1}(\theta|data)$, weighted by their relative posterior model probabilities, to obtain an appropriate ``averaged'' posterior distribution:
\begin{align}
\pi_{BMA}(\theta|data) =& \textrm{Pr}(M_{0}|data)\pi_{0}(\theta|data) +  \textrm{Pr}(M_{1}|data)\pi_{1}(\theta|data).
\label{eq:post_bma}
\end{align}
%

Alternatively, under the single ``mixture'' prior approach, one defines a single model with a single ``mixture'' prior, $\pi_{mix}(\theta)$, defined as a weighted combination of $\pi_{0}(\theta)$ and $\pi_{1}(\theta)$ such that:
\begin{equation}
 \pi_{mix}(\theta) = \textrm{Pr}(M_{0})\pi_{0}(\theta) +  \textrm{Pr}(M_{1})\pi_{1}(\theta).   \nonumber
\end{equation}
\noindent In this case, the corresponding posterior, $\pi_{mix}(\theta|data)$, can also be written as a weighted combination such that:
\begin{align}
    \pi_{mix}(\theta|data) &\propto  \Big(\textrm{Pr}(M_{0})\pi_{0}(\theta) +  \textrm{Pr}(M_{1})\pi_{1}(\theta)\Big)\times\textrm{Pr}(data|\theta).
\label{eq:post_mix}
\end{align}  
The $\pi_{mix}(\theta|data)$ posterior and the $\pi_{BMA}(\theta|data)$ posterior are in fact identical since:
\begin{align}    
      \pi_{mix}(\theta|data)  &\propto  \textrm{Pr}(M_{0})\pi_{0}(\theta)\textrm{Pr}(data|\theta) +  \textrm{Pr}(M_{1})\pi_{1}(\theta)\textrm{Pr}(data|\theta) \nonumber \\
      &\propto  m_{0}\textrm{Pr}(M_{0})\frac{\pi_{0}(\theta)\textrm{Pr}(data|\theta)}{m_{0}} +  m_{1}\textrm{Pr}(M_{1})\frac{\pi_{1}(\theta)\textrm{Pr}(data|\theta)}{m_{1}}, \nonumber  \\
         &\propto  \frac{m_{0}\textrm{Pr}(M_{0})\pi_{0}(\theta|data) +  m_{1}\textrm{Pr}(M_{1})\pi_{1}(\theta|data)}{m_{0}\textrm{Pr}(M_{0}) + m_{1}\textrm{Pr}(M_{1})} \nonumber  \\
     &\propto  m_{0}\textrm{Pr}(M_{0})\pi_{0}(\theta|data) +  m_{1}\textrm{Pr}(M_{1})\pi_{1}(\theta|data) \nonumber \\
              &\propto  \textrm{Pr}(M_{0}|data)\pi_{0}(\theta|data) +  \textrm{Pr}(M_{1}|data)\pi_{1}(\theta|data) \nonumber \\
              &\propto \pi_{BMA}(\theta|data),
\end{align}
\noindent where $m_{i} =\textrm{Pr}(data|M_{i})$, for $i=0,1$.  Thus, because densities are normalized to integrate to 1, we have that     $\pi_{mix}(\theta|data)=\pi_{BMA}(\theta|data)$.  
As such, point estimates and credible intervals obtained with the BMA approach will be identical to those obtained with the single ``mixture'' prior approach.   

This equality is not always acknowledged in the literature.  For example, in describing a mixture model for Bayesian meta-analysis, \citet{rover2019model} do note that the two approaches are equal (``the mixture prior effectively results in a model-averaging approach''), while  \citet{gronau2021primer} and \citet{bartovs2021bayesian} in proposing a similar approach for Bayesian meta-analysis, do not point out that the BMA approach is exactly equivalent to the single ``mixture'' prior approach.

Going forward we write simply $\pi(\theta|data)$ instead of using the $\pi_{mix}(\theta|data)$ or $\pi_{BMA}(\theta|data)$ notation, and use $\pi(\theta)$ to denote the ``mixed''/''averaged'' prior instead of $\pi_{mix}(\theta)$ or $\pi_{BMA}(\theta)$. In our example of the possibly biased coin, with the BMA approach, we obtain, from equation (\ref{eq:post_bma}):
\begin{align}
\pi(\theta|data) =& 0.523\times \delta_{0.5}(\theta) +  0.477\times f_{\textrm{Beta}}(\theta, 61,41). \nonumber
\end{align}
\noindent  Obtaining MCMC draws from this posterior is simple: With 1 million draws from $\pi_{0}(\theta|data)$ and another 1 million draws from $\pi_{1}(\theta|data)$, one can obtain 1 million draws from $\pi(\theta|data)$ by combining together approximately 523 thousand draws from  $\pi_{0}(\theta|data)$ with 477 thousand draws from $\pi_{1}(\theta|data)$.  See panel I of Figure \ref{fig:3by3} where the $\pi(\theta|data)$ function is plotted.  From these combined draws, we can then compute the posterior mean, $\hat{\theta}$, the posterior median, $\tilde{\theta}$, and the 95\% equal-tailed credible interval, $95\%\textrm{CrI}(\theta)$, as: $\hat{\theta}= 0.547$, $\tilde{\theta}= 0.500$, and $95\%\textrm{CrI}(\theta) = [0.500, 0.676]$.   Note that the 95\% equal-tailed credible interval includes 0.5.  \textcolor{black}{Also, note that due to the discontinuity in the posterior, this is a \textit{conservative} credible interval and is not actually equal-tailed: The [0.500, 0.676] interval includes 96.39\% of the posterior weight, since $\textrm{Pr}(\theta<0.500|data)=0.011$ and $\textrm{Pr}(\theta>0.676|data)=0.025$.}

Estimation based on $\pi(\theta|data)$ will be entirely congruent with testing based on the posterior model odds, $\textrm{PO}_{10}$, (or the Bayes factor, $\textrm{BF}_{10}$) since both estimation and testing are done using the very same (implied) prior model odds.  One could also calculate all of these numbers analytically.  For the posterior mean, we calculate:
\begin{align}
\hat{\theta} =& \quad \textrm{Pr}(M_{0}|data)\hat{\theta}_{0}  +  \textrm{Pr}(M_{1}|data)\hat{\theta}_{1} \nonumber \\
 =& \quad 0.523\times0.500 + 0.477\times0.598 \nonumber\\
 =&  \quad 0.547.
 \label{eq:postmean}
\end{align}
\noindent For the 95\% equal-tailed credible interval, we calculate $95\%\textrm{CrI}(\theta) = \Big[Q_{\theta|data}(0.025), Q_{\theta|data}(0.975)\Big]$, where:
\begin{align}
Q_{\theta|data}(q) =
\begin{cases}
& \quad Q_{\theta|data,M_{1}}\Big(\frac{q}{\textrm{Pr}(M_{1}|data)}\Big) ,  \quad \quad   \textrm{if}   \quad
\Big(\textrm{Pr}({\theta}<\theta_{0}|data, M_{1}) \textrm{Pr}(M_{1}|data)\Big) > q, \\
& \quad Q_{\theta|data,M_{1}}\Big(1-\frac{1-q}{\textrm{Pr}(M_{1}|data)}\Big) ,   \textrm{  if }   \Big(\textrm{Pr}({\theta}<\theta_{0}|data, M_{1}) \textrm{Pr}(M_{1}|data) + \textrm{Pr}(M_{0}|data)\Big) < q, \\
& \quad \theta_{0},  \quad \quad \quad \quad \quad \quad \quad \quad \quad \quad  \quad  \textrm{otherwise,}
\end{cases}
\label{eq:quantile}
\end{align}
\noindent where $Q_{Z}(q)$, is $q$-th quantile of $Z$, and $\theta_{0}=0.5$. Finally, for the posterior median, we calculate $\tilde{\theta} = Q_{\theta|data}(0.5)=0.500$.

\textcolor{black}{When $\theta_{0}$ is on the boundary of the 95\% credible interval, the $\Big[Q_{\theta|data}(0.025), Q_{\theta|data}(0.975)\Big]$ interval will necessarily be conservative (and not equal-tailed) in the sense that either $\textrm{Pr}(\theta<Q_{\theta|data}(0.025)|data)>0.025$ or $\textrm{Pr}(\theta>Q_{\theta|data}(0.975)|data)>0.025$ and therefore the interval will include more than 95\% of the posterior weight.  See \citet{campbell2022defining} for a discussion on this point.}

With the single ``mixture'' prior approach, from equation (\ref{eq:post_mix}), we obtain:
\begin{align}
    \pi(\theta|data) 
    &\propto  0.5\times f_{Binom}(X, N, \theta) \times \Big(\delta_{0.5}(\theta) + 1_{\{(0,1)\}}(\theta)\Big),
    \nonumber
\end{align}
 \noindent where $1_{\{(0,1)\}}(\theta)$ is an indicator function equal to 1 if $\theta\in(0,1)$ and equal to 0 otherwise, as in equation (\ref{eq:indicator}).  One might recognize this as  a version of the well-known ``spike-and-slab'' model (see \citet{van2021cautionary}) and Monte Carlo sampling directly from this posterior can be challenging when using popular MCMC software such as JAGS and stan.  An easy workaround is to introduce a latent parameter, $\omega$, such that the ``mixed prior'' is defined in a hierarchical way as follows:
\begin{align}
 \pi(\theta|\omega) &= (1-\omega)\pi_{0}(\theta) +  \omega\pi_{1}(\theta), \nonumber \\
   \omega &\sim \textrm{Bernoulli}(\textrm{Pr}(M_{1})). \nonumber
   \label{eq:hier}
\end{align}
\noindent This hierarchical strategy is often referred to as the ``product space method''; see \citet{carlin1995bayesian} and more recently \citet{lodewyckx2011tutorial}.   \textcolor{black}{See also the discussion about testing as mixture estimation in \citet{robert2016expected} and the discussion about unification via the spike-and-slab model in \citet{rouder2018bayesian}.}

Employing the ``product space method'' we obtain (using JAGS) 1 million draws from $\pi(\theta|data)$ and can calculate the posterior mean, $\hat{\theta}$, the posterior median, $\tilde{\theta}$, and the 95\% equal-tailed credible interval, $95\%\textrm{CrI}(\theta)$, as: $\hat{\theta}= 0.547$, $\tilde{\theta}= 0.500$, and $95\%\textrm{CrI}(\theta) = [0.500, 0.676]$.    These numbers are identical to those obtained using the BMA approach.  

 The ``product space method''  is also advantageous since $\hat{\omega}$, the posterior mean of $\omega$, is the posterior probability of $M_{1}$. As such, calculation of the posterior odds and of the Bayes factor is straightforward:
\begin{equation}
{\textrm{PO}}_{10} = \frac{\textrm{Pr}(M_{1}|data)}{\textrm{Pr}(M_{0}|data)} = \frac{\hat{\omega}}{1-\hat{\omega}}, \nonumber
\end{equation}
\noindent and:
%
\begin{equation}
{\textrm{BF}}_{10} = \left(\frac{\hat{\omega}}{1-\hat{\omega}}\right)/\left(\frac{\textrm{Pr}(M_{1})}{\textrm{Pr}(M_{0})}\right). \nonumber
\end{equation}

\vspace{0.1cm}

 For the possibly biased coin example, we obtain $\hat{\omega}=0.477$,  ${\textrm{PO}}_{10}= 0.913$, and ${\textrm{BF}}_{10} = 0.913$.  These numbers are equal the values calculated analytically in equations (\ref{eq:BF_coin}) and (\ref{eq:PO_coin}).




When the probability of $M_{0}$ given the data is non-negligible, ignoring the $M_{0}$ model ``paints an overly optimistic picture of what values $\theta$ is likely to have'' \citep{wagenmakers2020support}.  As such, there can be a substantial difference between the estimates obtained from the posterior under $M_{1}$ and estimates from the ``mixture''/``averaged'' posterior.  
In our example of the possibly biased coin, there is modest evidence in favour of $M_{0}$ with $\textrm{Pr}(M_{0}|data) = 0.523$, and as such there is a notable difference between $\hat{\theta}_{1}=0.598$ and $\hat{\theta}=0.547$, and between $95\%\textrm{CrI}(\theta)_{1}=[0.502,0.691]$ and  $95\%\textrm{CrI}(\theta)=[0.500,0.676]$.

\section{Observing a single coin flip}

\citet{wagenmakers2020principle} consider the coin flip example but with only a single flip (i.e., with $N=1$) which lands tails (i.e., with $X=0$).  In this simple, yet surprisingly interesting scenario, we have, from equation (\ref{eq:pr_data_m}):
\begin{equation}
    \textrm{Pr}(data|M_{0})= f_{Binom}(0, 1, 0.5)=0.500, \nonumber
\end{equation}
and
\begin{equation}
    \textrm{Pr}(data|M_{1})= \int_{t=0}^{t=1}f_{Binom}(0, 1, t)dt=0.500. \nonumber
\end{equation}
\noindent   The ratio of these two numbers is equal to the Bayes factor:
\begin{align}    
   \textrm{BF}_{10}    &= \frac{\textrm{Pr}(data|M_{1})}{\textrm{Pr}(data|M_{0})}  = \frac{0.500}{0.500} =  1,
   \nonumber
\end{align}
\noindent indicating that the single coin flip is ``perfectly uninformative'' \citep{wagenmakers2020principle}  with regards to determining whether the data support $M_{1}$ over $M_{0}$ or vice-versa. 

\citet{wagenmakers2020principle} do not explicitly define any prior model odds for this example, but do consider $\pi_{J}(\theta)$, a Beta(0.5,0.5) prior on $\theta$, for estimation.  Using this prior, a 95\% credible interval for $\theta$ is calculated as $95\%\textrm{CrI}(\theta)_{J}= [0.23, 1.00]$, and a 66\% credible interval for $\theta$ is calculated as $66\%\textrm{CrI}(\theta)_{J}= [0.70, 1.00]$.  \citet{wagenmakers2020principle} then explain their results as follows: ``It appears paradoxical that data can be perfectly uninformative for comparing [$M_{0}$ to $M_{1}$ ...], and at the same time provide reason to believe that $\theta > 0.5$ rather than $\theta < 0.5$.''  

\begin{figure}
    \centering
    \includegraphics[width=6.2in]{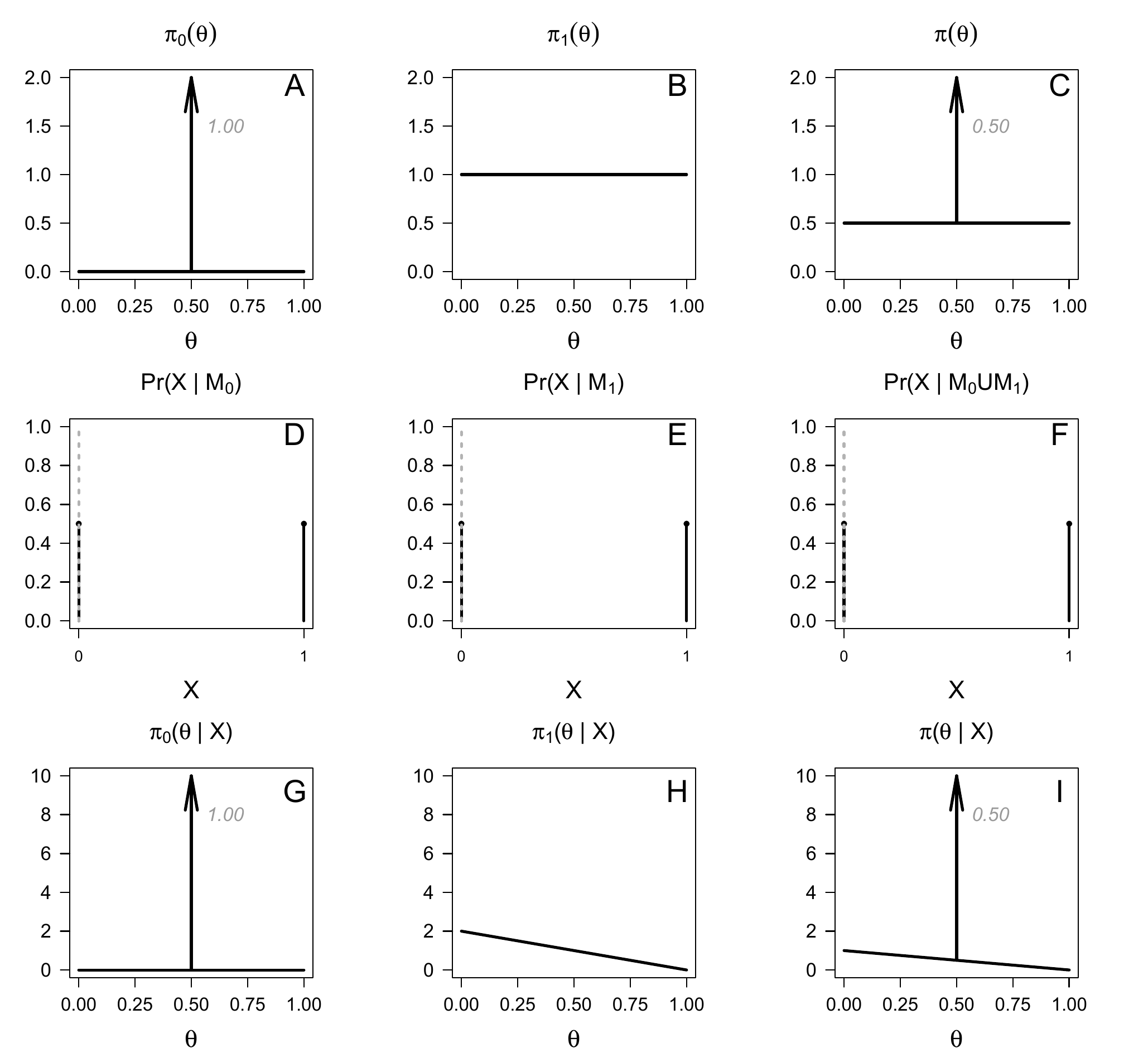}
    \caption{Observing a single coin flip ($X=0$ and $N=1$). The first row shows prior density functions (with the grey numbers next to the arrows corresponding to point mass of the ``spike''); middle row shows probability model for the data (with the grey vertical dashed lines corresponding to observed data of $X=60$); and the lower row shows posterior density functions (with the grey numbers next to the arrows corresponding to point mass of the ``spike'').  Left column corresponds to $M_{0}$, middle column corresponds to $M_{1}$, and right column corresponds to the ``mixed/averaged'' model.}
    \label{fig:single3by3}
\end{figure}

As was the case with $N=100$ coin flips, the ``paradox'' here is simply due to the fact that different implied prior model odds are being used for testing and for estimation.  If I have ``reason to believe'' that $\textrm{Pr}(\theta\in[0.70, 1.00]|data)=0.66$, then necessarily, I could not have believed, prior to flipping the coin, that either $M_{0}$ or $M_{1}$ were the true data generating mechanism.  This is because a belief that $\textrm{Pr}(\theta\in[0.70, 1.00]|data)=0.66$ necessarily implies a corresponding \textit{a priori} belief in the $\pi_{J}(\theta)$ prior.  

For someone who believes in the $\pi_{J}(\theta)$ prior, the above Bayes factor, $\textrm{BF}_{10}$, regardless of it's value, is meaningless since both $\textrm{Pr}(M_{0})=0$ and $\textrm{Pr}(M_{1})=0$.  The prior model odds will be ill-defined (i.e., $\textrm{PO}_{10}=$0:0) and therefore posterior model odds cannot be determined regardless of the observed data.  To be clear, the Bayes factor is only meaningful to someone if they are willing to consider that both models, $M_{0}$ and $M_{1}$, are \textit{a priori}, plausible.

Notably, if one assumes prior model odds of 1:1 (such that $\textrm{Pr}(M_{0})=\textrm{Pr}(M_{1}) = 0.5$), estimation based on the ``mixture''/``averaged'' posterior will be entirely congruent with the ``perfectly uninformative'' $\textrm{BF}_{10}=1$.  Indeed, we obtain $95\%\textrm{CrI}(\theta)= [0.025, 0.776]$ and a conservative $66\%\textrm{CrI}(\theta)= [0.188, 0.500]$,  which suggests that even a single flip is informative with regards to inferring $\theta$.  However, we also obtain:
\begin{align}
\textrm{Pr}(M_{0}|data)= \frac{0.5}{0.5 + 1 \times 0.5} = 0.500, \nonumber
\end{align}
\begin{align}
\textrm{Pr}(M_{1}|data)= \frac{1 \times 0.5 }{    0.5+1 \times 0.5} = 0.500, \quad 
\nonumber
\end{align}
and
\begin{align}
\textrm{PO}_{10} &=  \frac{0.5}{0.5} \times     1 = 1.000,
\nonumber
\end{align}
\noindent which suggests that, while we have gained some information about $\theta$, $M_{0}$ and $M_{1}$ remain equally likely to be the true model.  When framed in this manner, there is nothing paradoxical: Conclusions about $\theta$ and conclusions about $M_{0}$ vs. $M_{1}$ are entirely congruent since both are based on the same prior model odds.   Figure \ref{fig:single3by3} plots the ``mixture''/``averaged'' posterior where we see that, notably, exactly 50\% of posterior mass for $\theta$ is at exactly $\theta=0.500$.  Looking at the area under the curve, one can also determine that $\textrm{Pr}(\theta\ne0.5) = \textrm{Pr}(\theta<0.5) + \textrm{Pr}(\theta>0.5)=0.375 + 0.125 = 0.5$.

Using the ``mixture''/``averaged'' posterior for estimation (with the same prior model odds that are used for testing) will ensure that there is no ``incompatibility'' or ``mismatch'' between testing models and estimating $\theta$.  We emphasize that this applies universally, regardless of how $M_{0}$ and $M_{1}$ are defined.  Notably, some researchers have claimed that the ``incompatibility'' problem only occurs when considering point-null hypotheses (e.g.,  \citet{tendeiro2019review} write that: ``The problem is directly related to the use of the point null model'').  In the next section, we consider an example that does not involve a point null.

\section{With an interval null model}

Instead of questioning whether or not the coin is exactly fair, suppose one wishes to determine whether or not the coin is fair within some negligible margin.  Indeed, it could be argued that there is no such thing as a perfectly fair coin \citep{diaconis2007dynamical} and that bias is only consequential if it is of a certain non-negligible magnitude.

\begin{figure}
    \centering
    \includegraphics[width=6.2in]{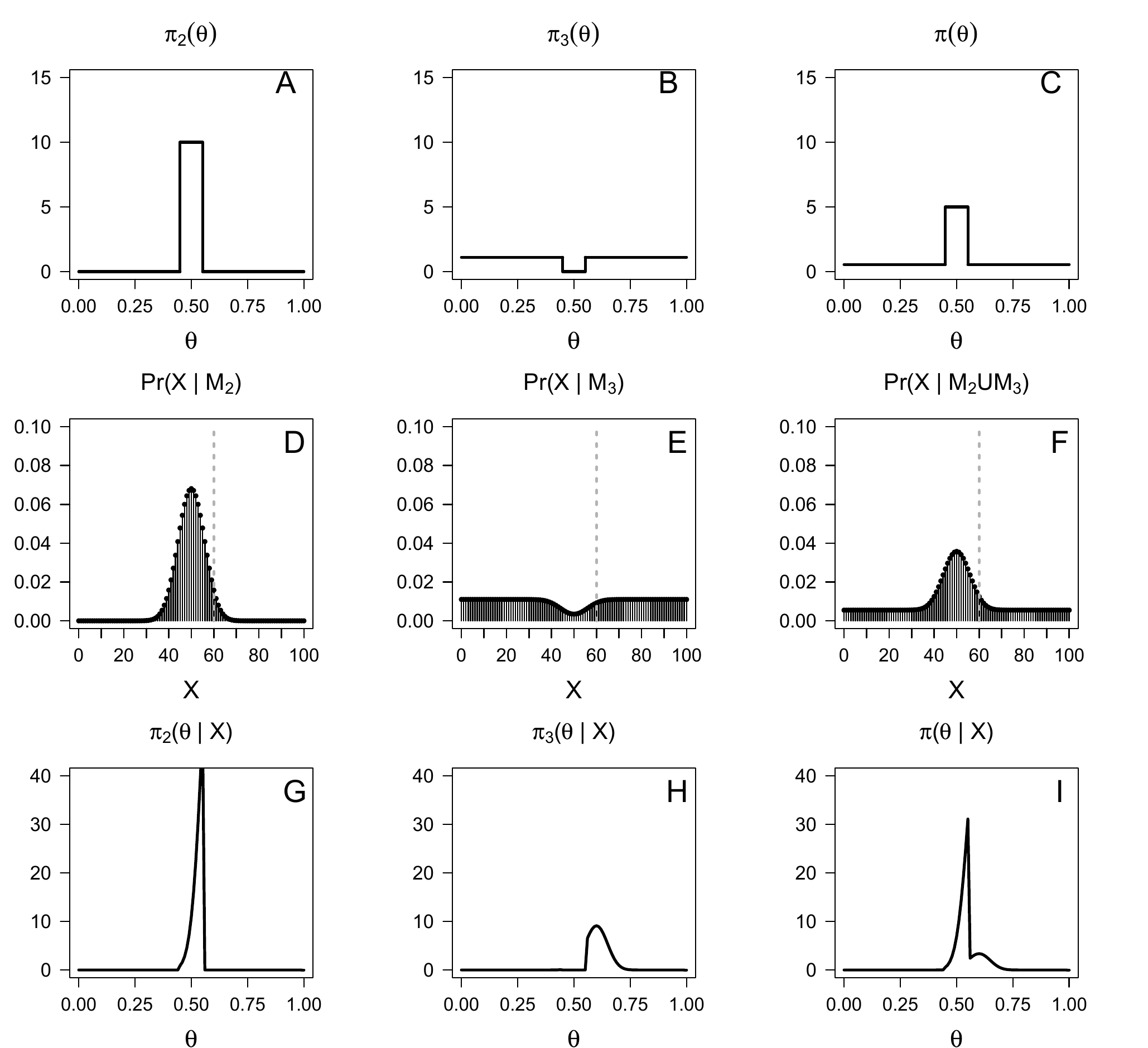}
    \caption{Observing a 100 coin flips ($X=60$ and $N=100$) with ``non-overlapping hypotheses''. The first row shows prior density functions; middle row shows probability model for the data (with the grey vertical dashed lines corresponding to observed data of $X=60$); and lower row shows posterior density functions.  Left column corresponds to $M_{2}$, middle column corresponds to $M_{3}$, and right column corresponds to the ``mixed/averaged'' model.}
    \label{fig:interval3by3}
\end{figure}

In other words, one might argue that $\theta$ will never be exactly equal to 0.5, and therefore consider ``non-overlapping hypotheses'' \citep{morey2011bayes} with priors for two competing models, $M_{2}$ and $M_{3}$, defined by partitioning a Uniform(0,1) density function as follows:
\begin{align}
 \pi_{2}(\theta)=
\begin{cases}
    10, & \text{if } \theta\in \Theta_{0} \\
    0,              & \text{otherwise,}
\end{cases}
\nonumber
\end{align}
and
\begin{align}
 \pi_{3}(\theta)=
\begin{cases}
    1/0.9,& \text{if } \theta\in \Theta_{1}  \\
    0,              & \text{otherwise.}
\end{cases}
\nonumber
\end{align}
\noindent where $\Theta_{0}=[0.45,0.55]$ and $\Theta_{1}=(0,0.45)\cup(0.55,1)$; see panels A and B of Figure \ref{fig:interval3by3}.  The corresponding posterior distributions (see panels G and H of Figure \ref{fig:interval3by3}) are:
\begin{align}
    \pi_{2}(\theta|data) &\propto f_{\textrm{Beta}}(\theta, X+1,N-X+1) \times 10\times 1_{\{\Theta_{0}\}}(\theta).
    \nonumber
\end{align}
%
%
\noindent and
\begin{align}
    \pi_{3}(\theta|data) &\propto f_{\textrm{Beta}}(\theta, X+1,N-X+1) \times 1.11\times 1_{\{\Theta_{1}\}}(\theta).
    \nonumber
\end{align}

Having observed $X=60$ ``heads'' out of $N=100$ coin flips, we calculate from equation (\ref{eq:pr_data_m}) (see panels D and E of Figure \ref{fig:interval3by3}):
\begin{equation}
    \textrm{Pr}(data|M_{2})= \int_{t=0.45}^{t=0.55}10\times f_{Binom}(60, 100, t)dt=0.016, 
    \nonumber
\end{equation}
\begin{equation}
    \textrm{Pr}(data|M_{3})= \frac{1}{0.9}\left(\int_{t=0}^{t=0.45}f_{Binom}(60, 100, t)dt + \int_{t=0.55}^{t=1}f_{Binom}(60, 100, t)dt\right) = 0.009.
    \nonumber
\end{equation}
\noindent and
\begin{align}    
   \textrm{BF}_{32}    &= \frac{\textrm{Pr}(data|M_{3})}{\textrm{Pr}(data|M_{2})}  = \frac{0.009}{0.016} =  0.584.
   \nonumber
\end{align}
%






For estimation, suppose one assumes a $\textrm{Uniform}(0,1)$ prior (i.e., the  $\pi_{1}(\theta)$ prior as defined in equation (\ref{eq:indicator})), so as to give equal prior weight to all values between 0 and 1.  From the corresponding $\pi_{1}(\theta|data)$ posterior (see equation (\ref{eq:post_m1})), one obtains a posterior mean of $\hat{\theta}_{1}=0.598$, a posterior median of $\tilde{\theta}_{1}=0.599$, and $95\%\textrm{CrI}(\theta)_{1}=[0.502,0.691]$.  One can also calculate $\textrm{Pr}_{1}(\theta \in \Theta_{1}|data)=0.840$.

This is therefore another example of the so-called ``mismatch'' between Bayesian testing and posterior estimation: The $ \textrm{BF}_{32}=0.584$ appears to support the null, while $95\%\textrm{CrI}(\theta)_{1}=[0.502,0.691]$ excludes 0.500.  As was the case in the previous examples with the point-null hypothesis, the ``mismatch'' here is due to different implied prior model odds being used for testing and for estimation.

As established by \citet{liao2021connecting}, the implied prior model odds that correspond to estimation with the Uniform(0,1) prior (i.e., under the $\pi_{1}(\theta|data)$ posterior) can be calculated as:
\begin{align}
   \textrm{IPMO}_{1} = \frac{\textrm{Pr}_{1}(\theta \in \Theta_{1})}{\textrm{Pr}_{1}(\theta \in \Theta_{0})} = 9:1,
\nonumber
\end{align}
 \noindent where 
\begin{align}
    \textrm{Pr}_{1}(\theta \in \Theta_{1}) =  \int_{\Theta_{1}}\pi_{1}(\theta)d\theta = \int_{t=0}^{t=0.45}1dt +\int_{t=0.55}^{t=1}1dt = 0.90.
    \nonumber
\end{align}
\noindent and
\begin{align}
    \textrm{Pr}_{1}(\theta \in \Theta_{0}) =  \int_{\Theta_{0}}\pi_{1}(\theta)d\theta = \int_{t=0.45}^{t=0.55}1dt = 0.10.
    \nonumber
\end{align}

\noindent Notably, by multiplying these implied prior model odds by the Bayes factor, one obtains posterior odds that favour $M_{3}$ and which are entirely congruent with the $\pi_{1}(\theta|data)$ posterior:
\begin{align}  
 \textrm{IPMO}_{1} \times \textrm{BF}_{32}  = 9 \times 0.584 = 5.255 = \frac{\textrm{Pr}_{1}(\theta \in \Theta_{1}|data)}{\textrm{Pr}_{1}(\theta \in \Theta_{0}|data)} = \frac{0.840}{0.160}.
 \nonumber
\end{align}

Alternatively, if one assumes equal prior model odds (i.e., assumes \textit{a priori} that $\textrm{Pr}(M_{2})=\textrm{Pr}(M_{3})=0.5$), then the posterior odds  favour $M_{2}$:
\begin{equation}
\textrm{PO}_{32} = \frac{\textrm{Pr}(M_{3})}{\textrm{Pr}(M_{2})} \times \textrm{BF}_{32}  =  \frac{0.5}{0.5} \times 0.584 = 0.584
\nonumber
\end{equation}

\noindent Estimation that is congruent with $\textrm{PO}_{32} =0.584$ can be done by invoking the ``mixed/averaged'' prior:
\begin{align}
{\pi}(\theta) &= \textrm{Pr}(M_{2})\pi_{2}(\theta) +\textrm{Pr}(M_{3})\pi_{3}(\theta), \nonumber \\
&= 0.5\times \frac{1}{0.1}\times1_{\{\Theta_{0}\}}(\theta) +  0.5\times \frac{1}{0.9}\times1_{\{\Theta_{1}\}}(\theta)
\nonumber
\end{align}
\noindent and the corresponding  ``mixed/averaged'' posterior:
\begin{align}
    \pi(\theta|data) &\propto \textrm{Pr}(M_{2}|data)\pi_{2}(\theta|data) + \textrm{Pr}(M_{3}|data)\pi_{3}(\theta|data) \nonumber \\
    &\propto 0.631  \times f_{\textrm{Beta}}(\theta, X+1,N-X+1)\times 10\times 1_{\{\Theta_{0}\}}(\theta) \quad + \nonumber  \\
    & \quad  \quad \quad \quad 0.369 \times f_{\textrm{Beta}}(\theta, X+1,N-X+1) \times 1.11\times 1_{\{\Theta_{1}\}}(\theta).
    \nonumber
\end{align}
\noindent  From the $\pi(\theta|data)$ posterior, one obtains a posterior mean of $\hat{\theta}=0.557$, a posterior median of $\tilde{\theta}=0.543$, and  $95\%\textrm{CrI}(\theta)=[0.479,0.673]$.  

To reiterate, testing and estimation will only be congruent if the (implied) prior model odds adopted for both are the same.  In this example, when estimation is based on a Uniform(0,1) prior for $\theta$, the implied prior model odds are 9:1, whereas when testing decisions about which model is most likely to be true are based on whether or not $\textrm{BF}_{32}>1$, the implied prior model odds are 1:1.  As such, it is no surprise that the two approaches provide different answers.

\section{A difference in personality types?}

As an example from applied research, consider \citet{heck2020review} who present the analysis of data concerning the difference between type D personality and non-type D personality for two outcome variables: BMI (body mass index) and NAS (negative affectivity score). See summary data in Table \ref{tab:heck} and note that these data were originally analyzed by \citet{lin2019negative}.  

 \citet{heck2020review} consider two competing models, a null model, $M_{0}$, which corresponds to the lack of an association between personality type and outcome, and an alternative model, $M_{1}$, which suggests the existence of an association.  Let $\theta$ be the parameter of interest, the standardized effect size.  The difference in personality types is associated with an increase of $\theta$ standard deviations in the outcome variable and a Normal model for the data is defined as:
\begin{align}
\textrm{Pr}(data|\theta, \sigma^{2}) =    \prod_{i=1}^{n} f_{Normal}\left (Y_{i},  \mu + X_{i}\times\frac{\sigma\theta}{2}, \sigma^{2}\right),
\nonumber
\end{align}
\noindent where $f_{Normal}()$ is the Normal probability density function,  $Y_{i}$ is the continuous outcome variable for the $i$-th observation (i.e., $Y_{i}$ is the $i$-th subject's BMI score in the first analysis, and the $i$-th subject's NAS score in the second analysis), and $X_{i}$ corresponds to the $i$-th subject's personality type such that $X_{i}=-1$ indicates a non-type D personality and $X_{i}=1$ indicates a type D personality.
 
 We have three parameters, $\mu$, $\theta$ and $\sigma^{2}$ for which we must define priors and we define these priors for both $M_{0}$ and $M_{1}$ as:
\begin{align}
\pi_{0}(\mu) = 1, \quad 
    \pi_{0}(\sigma^{2}) = \frac{1}{\sigma^{2}}, \quad \pi_{0}(\theta) = \delta_{0}(\theta), \quad
    \nonumber
\end{align}
and
\begin{align}
\pi_{1}(\mu) = 1, \quad 
    \pi_{1}(\sigma^{2}) = \frac{1}{\sigma^{2}}, \quad \pi_{1}(\theta) = f_\textrm{Cauchy}(\theta, 0.707), \quad
\nonumber
\end{align}
\noindent where, for both $M_{0}$ and $M_{1}$, noninformative Jeffreys priors are defined for $\sigma^{2}$, and for $\mu$; and $f_\textrm{Cauchy}(x,s) = ((s\pi)(1+(x/s)^{2}))^{-1}$ is the Cauchy probability density function evaluated at $x$, with scale parameter $s$. Specifying noninformative Jeffreys' priors is mathematically convenient for calculating the Bayes factor.  However, these improper priors cannot be accurately specified in popular MCMC sofware such as WinBugs, JAGS, and greta; see \citet{faulkenberry2018tutorial}.

For the BMI outcome, \citet{heck2020review} report a Bayes factor of $\textrm{BF}_{10}=1/4.2$ and sample from the posterior $\pi_{1}(\theta|data)$ (thereby implicitly assuming $M_{1}$ is correct) to obtain a posterior median of $\tilde{\theta}_{1}=0.05$, and  $95\%\textrm{CrI}(\theta)_{1}=[-0.35,0.45]$.  The posterior mean (not reported) is  $\hat{\theta}_{1}=0.05$.  

Assuming equal \textit{a priori} model probabilities, we can calculate from equation (\ref{eq:pr_m0_data}):
\begin{align}
\textrm{Pr}(M_{0}|data) \quad =& \quad \frac{\textrm{Pr}(M_{0})}{(BF_{10} \times \textrm{Pr}(M_{1}) + \textrm{Pr}(M_{0}))} \nonumber \\
=& \quad \frac{0.5}{0.5/4.2 + 0.5} \nonumber\\
=& \quad \frac{0.5}{0.5/4.2 + 0.5}\nonumber \\
=& \quad 0.808, \nonumber
\end{align}
indicating considerable evidence in favour of the null model.

Under $M_{0}$, we have $\theta=0$.  Therefore, we can obtain a Monte Carlo sample from the ``mixture''/``averaged'' posterior, $\pi(\theta|data)$, by simply combining 192 thousand draws for $\theta$ from $\pi_{1}(\theta|data)$ with 808 thousand zeros.  From this ``combined'' sample we calculate $\tilde{\theta}=0.00$, $\hat{\theta}=0.01$, and  $95\%\textrm{CrI}(\theta)=[-0.18, 0.28]$.  One can also obtain these estimates using equation (\ref{eq:quantile}) with $\theta_{0}=0$. These numbers are notably different than those based on the $\pi_{1}(\theta|data)$ posterior reported by \citet{heck2020review}. 

For the NAS outcome, \citet{heck2020review} report a Bayes factor of $\textrm{BF}_{10}=10^{20}$ and sample from the posterior $\pi_{1}(\theta|data)$ (implicitly assuming $M_{1}$ is correct) to obtain a posterior median of $\tilde{\theta}_{1}=-2.36$, and $95\%\textrm{CrI}(\theta)_{1}=[-2.85,-1.87]$.  The posterior mean (not reported) is  $\hat{\theta}_{1}=-2.36$. The evidence in favour of $M_{1}$ is overwhelming, $\textrm{Pr}(M_{1}|data)>0.99$, and therefore estimates from the ``mixture''/``averaged'' posterior (obtained using either the BMA approach or the ``mixture'' prior approach) are nearly identical: $\tilde{\theta}=-2.36$, $\hat{\theta}=-2.36$, and  $95\%\textrm{CrI}(\theta)=[-2.85,-1.87]$.

\begin{table}[h]
    \centering
    \begin{tabular}{|c|c|c|}
         \hline
        Outcome &  Non-Type D ($n=193$) & Type D ($n=23$)\\
         \hline
         BMI & $\bar{X}=26.0$ (SD=4.9) & $\bar{X}=25.7$ (SD=4.4)\\
         \hline
         NAS & $\bar{X}=5.5$ (SD=4.2) & $\bar{X}=15.4$ (SD=3.5) \\
              \hline
    \end{tabular}
    \caption{The sample mean, $\bar{X}$, and sample standard deviation, $SD$, for each of the two groups (type D personality and non-type D personality) and each of the two outcome variables (BMI (body mass index) and NAS (negative affectivity score)).  The data were originally analyzed by \citet{lin2019negative}.}
    \label{tab:heck}
\end{table}

\section{Flipping one million coins}

To better explain why we believe researchers should report estimates from the ``mixture''/``averaged'' posterior instead of estimates from $\pi_{1}(\theta|data)$ or other posteriors, we return one final time to the coin flip example, but this time with $N=10$ flips.  Let us explicitly assume that $M_{0}$ and $M_{1}$ are a \textit{priori} equally likely: $\textrm{Pr}(M_{0})=\textrm{Pr}(M_{1}) = 0.5$. 

Consider the following thought experiment.  By drawing from the ``mixture''/``averaged'' prior, $\pi(\theta)$, one obtains a coin for which $\theta=\theta^{[1]}$.  Flipping this coin $N=10$ times generates a total of $X=4$ ``heads'' and the posterior mean, under $M_{1}$, is:
\begin{align}
\hat{\theta}_{1}&=(X+\alpha)/(N+\alpha+\beta) \nonumber \\
 &= (4 + 1)/(10+1+1) \nonumber  \\
&= \frac{5}{12}=0.417.
\label{eq:x41}
\end{align}

The equal-tailed 95\% credible interval is $95\%\textrm{CrI}(\theta)_{1}=[0.167,0.692]$.  From equation (\ref{eq:postmean}), the posterior mean, taking into account the model uncertainty, is:
\begin{align}
\hat{\theta} &= \textrm{Pr}(M_{0}|data)\hat{\theta}_{0} + \textrm{Pr}(M_{1}|data)\hat{\theta}_{1} \nonumber\\
&=0.693\times0.5 + 0.307\times0.417\nonumber \\
&= 0.474, 
\label{eq:x4all}
\end{align}
and the equal-tailed 95\% credible interval is $95\%\textrm{CrI}(\theta)=[0.227,0.616]$. 
 
 Now suppose we obtain a second coin from the ``mixed''/``averaged'' prior for which $\theta=\theta^{[2]}$.  Flipping this second coin $N=10$ times generates a total of $X=8$ ``heads''.  Then the posterior mean, under $M_{1}$, is  $\hat{\theta}_{1}= \frac{9}{12}=0.750$, and the equal-tailed 95\% credible interval is $95\%\textrm{CrI}(\theta)_{1}=[0.482, 0.940]$.  \textcolor{black}{The posterior mean, taking into account the model uncertainty, is  $\hat{\theta}=0.669$, and the equal-tailed 95\% credible interval is $95\%\textrm{CrI}(\theta)=[0.500,0.930]$. Note that due to the discontinuity in the posterior, this is a \textit{conservative} 95\% credible interval that is not actually equal-tailed: The [0.500, 0.930] interval includes 95.29\% of the posterior weight, since $\textrm{Pr}(\theta<0.500|data)=0.022$ and $\textrm{Pr}(\theta>0.930|data)=0.025$.}

 We implemented this process by simulation again and again obtaining 1 million coins from the ``mixed''/``averaged'' prior (i.e., generating 1 million values $\theta^{[1]},\ldots,\theta^{[1000000]}$ from $\pi(\theta)$), and flipping each of these coins 10 times; see Table \ref{tab:millioncoins}.  In total, we obtained a total of 148,022 coins for which $X=4$; see Table \ref{tab:onetoten}.  
 
 Because this is a simulation, it is possible to reveal the underlying true values of $\theta^{[1]},\ldots,\theta^{[1000000]}$.  For the subset of coins for which $X=4$, the average of these values is 0.474.  For lack of better notation we write $\textrm{E}(\theta^{[j]}|X^{[j]}=4)= 0.474$. Notably, this equals $\hat{\theta}$ when $X=4$ (see equation (\ref{eq:x4all})), but does not equal $\hat{\theta}_{1}$ when $X=4$ (see equation (\ref{eq:x41})).  Nor does it equal $\frac{4}{10}$.  In this sense, both $\hat{\theta}_{1}$ and $\frac{X}{N}$ overestimate $\theta$, when $X=4$.  Table \ref{tab:onetoten} shows these comparisons for $X= 0,...,10$.
 
 Furthermore, amongst the 148,022 simulated $\theta^{[j]}$ values for which $X=4$, exactly 140,459 are within the interval of $95\%\textrm{CrI}(\theta)$=[0.227,0.616].  Notably, this proportion is 140,459/148,022 = 0.95.  In comparison, the proportion of the 148,022 values within the interval of $95\%\textrm{CrI}(\theta)_{1}$=[0.167, 0.692] is 0.98. Again, for lack of better notation, we write $\textrm{Pr}\Big(\theta^{[j]} \in 95\%\textrm{CrI}(\theta)_{1}\Big|X^{[j]}=4\Big)=0.98$,  and $\textrm{Pr}\Big(\theta^{[j]} \in 95\%\textrm{CrI}(\theta)\Big|X^{[j]}=4\Big)=0.95$.  In this sense, the $95\%\textrm{CrI}(\theta)_{1}$ interval is too wide, whereas the $95\%\textrm{CrI}(\theta)$ is appropriately wide.  Table \ref{tab:onetoten2} shows these comparisons for $X=0,...,10$.

\begin{table}[h]
    \centering
    \begin{tabular}{|c|c|c|c|c|c|c|c|}
    \hline
         $j$  & $\theta^{[j]}$ & $X^{[j]}$ & $N$ & $\hat{\theta}_{1}$ & $95\%\textrm{CrI}(\theta)_{1}$ & $\hat{\theta}$ & $95\%\textrm{CrI}(\theta)$\\
         \hline
         1 & 0.50 & 4 & 10 & 0.417  & [0.167,0.692] & 0.474   &   [0.227, 0.616]\\
         2 & 0.77 & 8 & 10 & 0.750 & [0.482, 0.940] & 0.669   &  [0.500,  0.930]\\
         3 & 0.91 & 8 & 10 & 0.750 & [0.482, 0.940] & 0.669   &  [0.500,  0.930]\\
         4 & 0.11 & 2 & 10 & 0.250 & [0.060, 0.518] & 0.331&   [0.070, 0.500]\\
         5 & 0.87 & 9 & 10 & 0.833 & [0.587, 0.977] & 0.801&   [0.500, 0.976]\\
         6 & 0.50 & 7 & 10 & 0.667 & [0.390, 0.891] & 0.573&   [0.445, 0.860]\\
         7 & 0.50 & 5 & 10 & 0.500 & [0.234, 0.766] & 0.500&   [0.312, 0.688]\\
         8 & 0.31 & 1 & 10 & 0.167 & [0.023, 0.413] & 0.199&   [0.024, 0.500]\\  
         $\vdots$ & $\vdots$ & $\vdots$ & $\vdots$ & $\vdots$ & $\vdots$  & $\vdots$  & $\vdots$ \\
         $1,000,000$ & 0.60 & 6 & 10 & 0.583 & [0.308, 0.833] & 0.526  &   [0.384, 0.773]\\  
         \hline
    \end{tabular}
    \caption{One million coins: For the $j$-th coin we observe $X^{[j]}$ ``heads'' out of $N=10$ flips, we calculate the posterior estimates.   \textcolor{black}{Whenever $0.500$ is on the boundary of the $95\%\textrm{CrI}(\theta)$ interval, the interval will be conservative.  Specifically, for $X=1$ and $X=9$, 96.97\% of the posterior weight will be within the interval; for $X=2$ and $X=8$, 95.29\% of the posterior weight will be within the interval.}}
    \label{tab:millioncoins}
\end{table}

\begin{table}[ht]
\centering
\begin{tabular}{|r|rcccc|}
  \hline
$x$ & \# of coins & $\frac{x}{N}$ & $\hat{\theta}_{1}$ &$\hat{\theta}$ & $\textrm{E}(\theta^{[j]}|X^{[j]}=x)$ \\ 
  \hline
0 & 45,864 & 0.0 & 0.083 & 0.088 & 0.087 \\ 
  1 & 50,180 & 0.1 & 0.167 & 0.199 & 0.198 \\ 
  2 & 67,281 & 0.2 & 0.250 & 0.331 & 0.331 \\ 
  3 & 104,309 & 0.3 & 0.333 & 0.427 & 0.427 \\ 
  4 & 148,022 & 0.4 & 0.417 & 0.474 & 0.474 \\ 
  5 & 168,651 & 0.5 & 0.500 & 0.500 & 0.500 \\ 
  6 & 148,502 & 0.6 & 0.583 & 0.526 & 0.526 \\ 
  7 & 103,961 & 0.7 & 0.667 & 0.573 & 0.573 \\ 
  8 & 67,353 & 0.8 & 0.750 & 0.669 & 0.670 \\ 
  9 & 49,939 & 0.9 & 0.833 & 0.801 & 0.801 \\ 
  10 & 45,938 & 1.0 & 0.917 & 0.912 & 0.913 \\ 
   \hline
\end{tabular}
\caption{Simulation study results: Average values obtained from the subset set of ``coins'' for which  $X=x$ ``heads''.  For example, among the subset of $104,309$ coins that produced $X=3$, the average of $\theta^{[j]}$ values
is $0.427$, a numerical approximation to $\textrm{E}(\theta^{[j]}|X^{[j]}=3)$. Whereas the posterior mean estimate (taking into account the model uncertainty) whenever $X^{[j]}=3$ is $\hat{\theta}=0 .427$.}
\label{tab:onetoten}
\end{table}

\begin{table}[ht]
\centering

\begin{tabular}{|r|rcc|}
  \hline
$x$ & \# of coins  & $\textrm{Pr}\Big(\theta^{[j]} \in 95\%\textrm{CrI}(\theta)_{1}|X^{[j]}=x\Big)$  & $\textrm{Pr}\Big(\theta^{[j]} \in 95\%\textrm{CrI}(\theta)|X^{[j]}=x\Big)$ \\ 
  \hline
0 & 45,864 & 0.940 & 0.950 \\ 
  1 & 50,180 & 0.858 & 0.968 \\ 
  2 & 67,281 & 0.967 & 0.953 \\ 
  3 & 104,309 & 0.978 & 0.951 \\ 
  4 & 148,022 & 0.984 & 0.949 \\ 
  5 & 168,651 & 0.987 & 0.951 \\ 
  6 & 148,502 & 0.985 & 0.950 \\ 
  7 & 103,961 & 0.978 & 0.950 \\ 
  8 & 67,353 & 0.967 & 0.955 \\ 
  9 & 49,939 & 0.860 & 0.970 \\ 
  10 & 45,938 & 0.938 & 0.949 \\ 
   \hline
\end{tabular}
\caption{Simulation study results: Proportion of coins (from the subset for which $X=x$ ``heads'') within each credible interval. For example, amongst the 148,502 simulated $\theta^{[j]}$ values for which $X=6$, exactly 95.0\% are within the interval of $95\%\textrm{CrI}(\theta)$=[0.384, 0.773].  \textcolor{black}{Note that due to the discontinuity in the posterior, for each of $X=1$, $X=2$, $X=8$, and $X=9$, the $95\%\textrm{CrI}(\theta)$ interval is \textit{conservative} and therefore we expect to see that more than 95.0\% of simulated $\theta^{[j]}$ values are within the interval; see caption of Table \ref{tab:millioncoins} for specific values.}}
\label{tab:onetoten2}
\end{table}

The behaviour exhibited in Tables \ref{tab:onetoten} and \ref{tab:onetoten2} is completely general mathematically.  Generically, let $(\theta ^{[j]},X ^{[j]})$ be a random draw as is used to construct Table \ref{tab:millioncoins}, row by row, e.g., $\theta ^{[j]}$ is drawn from $\pi(\theta)$, and then $X ^{[j]}$ is drawn from ``the model'' with $\theta$ set equal to $\theta ^{[j]}$.   
Thus the average of $\theta^{[j]}$ when $X=x$ (the last column of Table \ref{tab:onetoten}) is simply (a Monte Carlo approximation to) $\textrm{E}(\theta ^{[j]}|X ^{[j]}=x)$, for each $x$.  But since  $(\theta ^{[j]},X^{[j]})$ is nothing more than a joint draw from the amalgamation of the prior and statistical model, 
$\textrm{E}(\theta ^{[j]}|X ^{[j]}=x)$ and $\hat{\theta}$ when $X=x$, are one and the same.  By definition then, if we repeatedly sample parameter and data pairs this way, the average of the parameter values amongst draws yielding a specific data value {\em is} the posterior mean of the parameter given that data value.  The same calibration idea applies to credible intervals so that $\textrm{Pr}\Big(\theta^{[j]} \in 95\%\textrm{CrI}(\theta)\Big|X^{[j]}=x\Big)\ge0.95$. Note that in the absence of a point mass (i.e., without the ``spike''), this will be a strict equality \citep{campbell2022defining}.

This strongly supports the argument that when considering model uncertainty via a Bayes factor, $\hat{\theta}$ is the appropriate estimate to report and $95\%\textrm{CrI}(\theta)$ is the appropriate credible interval.  While these sorts of calibration properties of Bayesian estimators with respect to  repeated sampling of data under {\em different}  parameter values feature sporadically in the literature (e.g., \citet{rubin1986efficiently}), they do not seem to be widely known. They do however, underpin the scheme of \citet{cook2006validation} to validate Markov Chain Monte Carlo approximations to posterior quantities and they are also discussed more recently in \citet{gustafson2009interval}.  Importantly, we note that this sense of calibration arises only when the prior truly corresponds to the data-generating process.   

\section{Conclusion}

In conclusion, there is no mismatch between the Bayes factor (or the posterior odds) and Bayesian posterior estimation with an appropriately defined prior and posterior.  Our recommendation is therefore very simple.  If one reports a Bayes factor comparing $M_{0}$ and $M_{1}$, then one should also report posterior estimates based on the ``mixed/``averaged'' posterior, $\pi(\theta|data)$, with prior model odds appropriately specified.  Researchers should refrain from reporting posterior estimates based on $\pi_{1}(\theta|data)$ (which implicitly assumes $M_{1}$ is correct) or estimates based on other posteriors for which the Bayes factor is meaningless.  We see no reason why disregarding $M_{0}$ ``for the purpose of parameter estimation''  \citep{wagenmakers2020overwhelming} is advisable.   

\textcolor{black}{On a final note, in order to keep things as accessible as possible, all of the examples considered in this paper were univariate inference problems involving only two models.  Our conclusions and recommendations, however, also apply to more complex problems.  On this point, we refer interested readers to \citet{rossell2017nonlocal}, who consider strategies for how to use the same prior for both estimation and model selection in high-dimensional settings to achieve coherence between estimation and testing, and to \citet{lavine1999bayes} and \citet{oelrich2020bayesian}, who discuss Bayes factor based model selection in scenarios involving more than two models.}

\bibliography{truthinscience}

@article{wagenmakers2020support,
  title={The support interval},
  author={Wagenmakers, Eric-Jan and Gronau, Quentin F and Dablander, Fabian and Etz, Alexander},
  journal={Erkenn},
  volume={87},
  pages={589--601},
  year={2022},

  publisher={Springer}
}

@article{rouder2018bayesian,
  title={Bayesian inference for psychology, part {I}{V}: {P}arameter estimation and {B}ayes factors},
  author={Rouder, Jeffrey N and Haaf, Julia M and Vandekerckhove, Joachim},
  journal={Psychonomic {B}ulletin \& {R}eview},
  volume={25},
  number={1},
  pages={102--113},
  year={2018},
  publisher={Springer}
}

@misc{wagenmakers2020overwhelming,
  title={Overwhelming Evidence for Vaccine Efficacy in the {P}fizer Trial: An Interim {B}ayesian Analysis},
  author={Wagenmakers, Eric-Jan and Gronau, Quentin Frederik},
  year={2020},
  howpublished={PsyArXiv: 10.31234/osf.io/fs562}
}

@article{rossell2017nonlocal,
  title={Nonlocal priors for high-dimensional estimation},
  author={Rossell, David and Telesca, Donatello},
  journal={Journal of the {A}merican {S}tatistical {A}ssociation},
  volume={112},
  number={517},
  pages={254--265},
  year={2017},
  publisher={Taylor \& Francis}
}

@incollection{wagenmakers2020principle,
  title={The principle of predictive irrelevance or why intervals should not be used for model comparison featuring a point null hypothesis},
  author={Wagenmakers, Eric-Jan and Lee, Michael D and Rouder, Jeffrey N and Morey, Richard D},
  booktitle={The {T}heory of {S}tatistics in {P}sychology},
  pages={111--129},
  year={2020},
  publisher={Springer},
  address={Cham, Switzerland},
  editor={Craig W. Gruber},
}

@article{van2021advantages,
  title={Advantages masquerading as' issues' in {B}ayesian hypothesis testing: {A} commentary on {T}endeiro and {K}iers (2019).},
  author={van Ravenzwaaij, Don and Wagenmakers, Eric-Jan},
  journal={Psychological {M}ethods},
  year={2021},
  volume={27},
  number={3}, 
  pages={451--465}
}

@article{robert2016expected,
  title={The expected demise of the {B}ayes factor},
  author={Robert, Christian P},
  journal={Journal of Mathematical Psychology},
  volume={72},
  pages={33--37},
  year={2016},
  publisher={Elsevier}
}

@article{wagenmakers2018bayesian,
  title={Bayesian inference for psychology. {P}art {I}: {T}heoretical advantages and practical ramifications},
  author={Wagenmakers, Eric-Jan and Marsman, Maarten and Jamil, Tahira and Ly, Alexander and Verhagen, Josine and Love, Jonathon and Selker, Ravi and Gronau, Quentin F and {\v{S}}m{\'\i}ra, Martin and Epskamp, Sacha and others},
  journal={Psychonomic {B}ulletin \& {R}eview},
  volume={25},
  number={1},
  pages={35--57},
  year={2018},
  publisher={Springer}
}

@article{tendeiro2019review,
  title={A review of issues about null hypothesis Bayesian testing.},
  author={Tendeiro, Jorge N and Kiers, Henk AL},
  journal={Psychological {M}ethods},
  volume={24},
  number={6},
  pages={774--795},
  year={2019},
  publisher={American Psychological Association}
}

@article{morey2011bayes,
  title={Bayes factor approaches for testing interval null hypotheses.},
  author={Morey, Richard D and Rouder, Jeffrey N},
  journal={Psychological {M}ethods},
  volume={16},
  number={4},
  pages={406--419},
  year={2011},
  publisher={American Psychological Association}
}

@article{liao2021connecting,
  title={Connecting and contrasting the Bayes factor and a modified ROPE procedure for testing interval null hypotheses},
  author={Liao, JG and Midya, Vishal and Berg, Arthur},
  journal={The {A}merican {S}tatistician},
  volume={75},
  number={3},
  pages={256--264},
  year={2021},
  publisher={Taylor \& Francis}
}

@article{diaconis2007dynamical,
  title={Dynamical bias in the coin toss},
  author={Diaconis, Persi and Holmes, Susan and Montgomery, Richard},
  journal={SIAM review},
  volume={49},
  number={2},
  pages={211--235},
  year={2007},
  publisher={SIAM}
}

@article{puga2015bayesian,
  title={Bayesian statistics: today's predictions are tomorrow's priors},
  author={Puga, Jorge Lopez and Krzywinski, Martin and Altman, Naomi},
  journal={Nature Methods},
  volume={12},
  number={5},
  pages={377--379},
  year={2015},
  publisher={Nature Publishing Group}
}

@book{o2004kendall,
  title={Kendall's advanced theory of statistics, volume 2B: {B}ayesian inference},
  author={O'Hagan, Anthony and Forster, Jonathan J},
  edition={2},
  year={2004},
  address={London},
  publisher={Arnold}
}

@article{lavine1999bayes,
  title={Bayes factors: {W}hat they are and what they are not},
  author={Lavine, Michael and Schervish, Mark J},
  journal={The {A}merican {S}tatistician},
  volume={53},
  number={2},
  pages={119--122},
  year={1999},
  publisher={Taylor \& Francis}
}

@article{puga2015bayes,
  title={Bayes' theorem: {I}ncorporate new evidence to update prior information},
  author={Puga, Jorge Lopez and Krzywinski, Martin and Altman, Naomi},
  journal={Nature {M}ethods},
  volume={12},
  number={4},
  pages={277--279},
  year={2015},
  publisher={Nature Publishing Group}
}

@article{oelrich2020bayesian,
  title={When are {B}ayesian model probabilities overconfident?},
  author={Oelrich, Oscar and Ding, Shutong and Magnusson, M{\aa}ns and Vehtari, Aki and Villani, Mattias},
  journal={arXiv preprint arXiv:2003.04026},
  year={2020}
}

@article{campbell2022defining,
  title={Defining a credible interval is not always possible with ``point-null'' priors: {A} lesser-known consequence of the {J}effreys-{L}indley paradox},
  author={Campbell, Harlan and Gustafson, Paul},
  howpublished={arxiv},
  doi={https://doi.org/10.48550/arXiv.2210.00029},
  year={2022},
  journal = {arXiv preprint  arXiv:2210.00029},
   eprint = {math.ST/2210.00029v1}
}

@article{rubin1986efficiently,
  title={Efficiently simulating the coverage properties of interval estimates},
  author={Rubin, Donald B and Schenker, Nathaniel},
  journal={Journal of the {R}oyal {S}tatistical {S}ociety: {S}eries {C} ({A}pplied {S}tatistics)},
  volume={35},
  number={2},
  pages={159--167},
  year={1986},
  publisher={Wiley Online Library}
}

@article{lodewyckx2011tutorial,
  title={A tutorial on {B}ayes factor estimation with the product space method},
  author={Lodewyckx, Tom and Kim, Woojae and Lee, Michael D and Tuerlinckx, Francis and Kuppens, Peter and Wagenmakers, Eric-Jan},
  journal={Journal of {M}athematical {P}sychology},
  volume={55},
  number={5},
  pages={331--347},
  year={2011},
  publisher={Elsevier}
}

@article{rover2019model,
  title={Model averaging for robust extrapolation in evidence synthesis},
  author={R{\"o}ver, Christian and Wandel, Simon and Friede, Tim},
  journal={Statistics in {M}edicine},
  volume={38},
  number={4},
  pages={674--694},
  year={2019},
  publisher={Wiley Online Library}
}

@article{faulkenberry2018tutorial,
  title={A tutorial on generalizing the default {B}ayesian $t$-test via posterior sampling and encompassing priors},
  author={Faulkenberry, Thomas J},
  journal={arXiv preprint arXiv:1812.03092},
  year={2018}
}

@article{bartovs2021bayesian,
  title={Bayesian model-averaged meta-analysis in medicine},
  author={Barto{\v{s}}, Franti{\v{s}}ek and Gronau, Quentin F and Timmers, Bram and Otte, Willem M and Ly, Alexander and Wagenmakers, Eric-Jan},
  journal={Statistics in {M}edicine},
  volume={40},
  number={30},
  pages={6743--6761},
  year={2021},
  publisher={Wiley Online Library}
}

@article{gronau2021primer,
  title={A primer on {B}ayesian model-averaged meta-analysis},
  author={Gronau, Quentin F and Heck, Daniel W and Berkhout, Sophie W and Haaf, Julia M and Wagenmakers, Eric-Jan},
  journal={Advances in Methods and Practices in Psychological Science},
  volume={4},
  number={3},
  doi={25152459211031256},
  pages={doi:10.1177/25152459211031256},
  year={2021},
  publisher={Sage Publications Sage CA: Los Angeles, CA}
}

@article{carlin1995bayesian,
  title={Bayesian model choice via {M}arkov chain {M}onte {C}arlo methods},
  author={Carlin, Bradley P and Chib, Siddhartha},
  journal={Journal of the {R}oyal {S}tatistical {S}ociety: {S}eries {B} ({M}ethodological)},
  volume={57},
  number={3},
  pages={473--484},
  year={1995},
  publisher={Wiley Online Library}
}

@article{lin2019negative,
  title={Negative affectivity and social inhibition are associated with increased cardiac readmission in patients with heart failure: {A} preliminary observation study},
  author={Lin, Tin-Kwang and You, Kai-Xun and Hsu, Chiu-Tien and Li, Yi-Da and Lin, Chin-Lon and Weng, Chia-Ying and Koo, Malcolm},
  journal={Plos {O}ne},
  volume={14},
  number={4},
  pages={e0215726},
  year={2019},
  publisher={Public Library of Science San Francisco, CA USA}
}

@article{gelman2002you,
  title={You can load a die, but you can't bias a coin},
  author={Gelman, Andrew and Nolan, Deborah},
  journal={The {A}merican {S}tatistician},
  volume={56},
  number={4},
  pages={308--311},
  year={2002},
  publisher={Taylor \& Francis}
}

@article{van2021cautionary,
  title={A cautionary note on estimating effect size},
  author={van den Bergh, Don and Haaf, Julia M and Ly, Alexander and Rouder, Jeffrey N and Wagenmakers, Eric-Jan},
  journal={Advances in Methods and Practices in Psychological Science},
  volume={4},
  number={1},
  pages={doi:10.1177/2515245921992035},
  year={2021},
  publisher={Sage Publications Sage CA: Los Angeles, CA}
}

@article{cook2006validation,
title={Validation of software for {B}ayesian models using posterior quantiles},
author={Cook, Samantha R and Gelman, Andrew and Rubin, Donald B},
journal={Journal of {C}omputational and {G}raphical {S}tatistics},
volume={15},
number={3},
pages={675--692},
year={2006},
publisher={Taylor \& Francis}
}

@article{gustafson2009interval,
 title={Interval estimation for messy observational data},
 author={Gustafson, Paul and Greenland, Sander},
 journal={Statistical {S}cience},
 volume={24},
 number={3},
 pages={328--342},
 year={2009},
 publisher={Institute of Mathematical Statistics}
}

@article{kelter2022evidence,
  title={The evidence interval and the {B}ayesian evidence value: {O}n a unified theory for {B}ayesian hypothesis testing and interval estimation},
  author={Kelter, Riko},
  journal={British {J}ournal of {M}athematical and {S}tatistical {P}sychology},
   volume={75}, 
   pages={550--592},
  year={2022},
  publisher={Wiley Online Library}
}

@article{heck2020review,
  title={A review of applications of the {B}ayes factor in psychological research.},
  author={Heck, Daniel W and Boehm, Udo and B{\"o}ing-Messing, Florian and B{\"u}rkner, Paul-Christian and Derks, Koen and Dienes, Zoltan and Fu, Qianrao and Gu, Xin and Karimova, Diana and Kiers, Henk AL and others},
  journal={Psychological Methods},
  year={2022},
  number={Advance online publication},
  publisher={American Psychological Association}
}

\end{document}